\documentclass[10pt,fleqn]{article}
\usepackage{bbm}
\usepackage{amsfonts}
\usepackage{amssymb}
\usepackage{latexsym,amsmath}
\usepackage{graphicx}

\date{}
\begin{document}
\title{Bounds on the spectral radius of nonnegative matrices and applications in graph spectra}
\author{Shu-Yu Cui $^a$, Gui-Xian Tian $^b$\footnote{Corresponding author.  E-mail: gxtian@zjnu.cn or guixiantian@163.com (G.-X.
Tian)}\\%EndAName
{\small{\it $^a$Xingzhi College, Zhejiang Normal University, Jinhua,
Zhejiang, 321004, P.R. China}} \\
{\small{\it $^b$College of Mathematics, Physics and Information Engineering,}}\\
{\small{\it Zhejiang Normal University, Jinhua, Zhejiang, 321004,
P.R. China}} }\maketitle

\begin{abstract}

In this paper, we give upper and lower bounds for the spectral
radius of a nonnegative irreducible matrix and characterize the
equality cases. These bounds theoretically improve and generalize
some known results of Duan et al.[X. Duan, B. Zhou, Sharp bounds on
the spectral radius of a nonnegative matrix, Linear Algebra Appl.
(2013), http://dx.doi.org/10.1016/j.laa.2013.08.026]. Finally,
applying these bounds to various matrices associated with a graph,
we obtain some new upper and lower bounds on various spectral
radiuses of graphs, which generalize and improve some known results.


\end{abstract}

\section*{1. Introduction}

We only consider simple graphs which have no loops and multiple
edges. Let $G=(V,E)$ be a simple graph with vertex set
$V=\{v_1,v_2,\ldots,v_n\}$ and edge set $E$. For any two vertices
$v_i,v_j\in V$, we write $i\sim j$ if $v_i$ and $v_j$ are adjacent.
For any vertex $v_i\in V$, denote the \emph{degree} of $v_i$ by
$d_i$.

The matrix $A(G)=(a_{ij})$ is called the \emph{adjacency matrix} of
$G$, where $a_{ij}=1$ if $i\sim j$ and $0$ otherwise. Let
$D(G)=\text{diag}(d_1,d_2,\ldots,d_n)$ be the diagonal matrix of
vertex degrees. Then $Q(G)=D(G)+A(G)$ is called the \emph{signless
Laplacian matrix} of $G$. For more review about the signless
Laplacian matrix of $G$, readers may refer to
\cite{Cvetkovic2008,CvetkovicI,CvetkovicII,CvetkovicIII} and the
references therein.

If $G$ is connected, then the \emph{distance matrix} of $G$ is the
$n\times n$ matrix $\mathbb{D}(G) = (d_{ij})$, where $d_{ij}$ is the
distance between two vertices $v_i$ and $v_j$ in $G$. For any vertex
$v_i\in V$, denote the \emph{transmission} of $v_i$ by
$\mathbb{D}_i$, i.e., the the sum of distances between $v_i$ and
other vertices of $G$. Let
$\mathbb{T}\mathbbm{r}(G)=\text{diag}(\mathbb{D}_1,\ldots,\mathbb{D}_n)$
be the transmission diagonal matrix. The \emph{distance signless
Laplacian matrix} of $G$ is the $n\times n$ matrix
$\mathbb{D}\mathbb{Q}(G)=\mathbb{T}\mathbbm{r}(G)+\mathbb{D}(G)$.
For more review about the distance signless Laplacian matrix of $G$,
see \cite{Aouchiche2013}.

Let A=$(a_{ij})$ be an $n\times n$ nonnegative matrix. The spectral
radius of $A$, denoted by $\rho(A)$, is the largest modulus of an
eigenvalue of $A$. Moreover, if $A$ is symmetric, then $\rho(A)$ is
just the largest eigenvalue of $A$. For more review about
nonnegative matrix, see \cite{Berman1994,Minc1988}.

In \cite{Chen2013}, Chen et al. obtained some upper bounds for the
spectral radius of a nonnegative irreducible matrix with diagonal
entries 0 and characterized the equality cases. Using the same
technique, Duan and Zhou\cite{Duan2013} got some upper and lower
bounds for the spectral radius of a nonnegative matrix and
characterized the equality cases whenever the matrix is irreducible
in the general case. Applying these bounds to the matrices
associated with a graph as mentioned above, they also obtained some
upper and lower bounds on various spectral radiuses of graphs. For
example,

Let $d_1\geq d_2\geq \cdots \geq d_n$ be the degree sequence of $G$.
Then, for $1\leq i\leq n$,

\begin{equation}\label{1}
\rho (A(G)) \le  {\frac{{d_i - 1 + \sqrt {(d_i + 1)^2 +
4\sum\nolimits_{k = 1}^{i - 1} {(d_k - d_i )} } }}{2}}
\cite{Duan2013,Liu2013},
\end{equation}

\begin{equation}\label{2}
\rho (Q(G)) \le {\frac{{d_1  + 2d_i  - 1 + \sqrt {(2d_i  - d_1  +
1)^2  + 8\sum\nolimits_{k = 1}^{i - 1} {(d_k  - d_i )} } }}{2}}
\cite{Cui2013,Duan2013}.
\end{equation}

Let $\mathbb{D}_1\geq \mathbb{D}_2\geq \cdots \geq \mathbb{D}_n$ be
the transmission sequence of $G$ and $d$ be the diameter of $G$.
Then, for $1\leq i\leq n$,

\begin{equation}\label{3}
\rho (\mathbb{D}(G)) \le  {\frac{{\mathbb{D}_i - d + \sqrt
{(\mathbb{D}_i + d)^2 + 4d\sum\nolimits_{k = 1}^{i - 1}
{(\mathbb{D}_k - \mathbb{D}_i )} } }}{2}} \cite{Chen2013,Duan2013},
\end{equation}

\begin{equation}\label{4}
\rho (\mathbb{D}\mathbb{Q}(G)) \le {\frac{{\mathbb{D}_1  +
2\mathbb{D}_i  - d + \sqrt {(2\mathbb{D}_i  - \mathbb{D}_1  + d)^2 +
8d\sum\nolimits_{k = 1}^{i - 1} {(\mathbb{D}_k  - \mathbb{D}_i )} }
}}{2}} \cite{Duan2013}
\end{equation}
and

\begin{equation}\label{5}
\rho (\mathbb{D}(G)) \geq {\frac{{\mathbb{D}_n - 1 + \sqrt
{(\mathbb{D}_n+ 1)^2 + 4\sum\nolimits_{k = 1}^{n - 1} {(\mathbb{D}_k
- \mathbb{D}_n )} } }}{2}} \cite{Duan2013},
\end{equation}

\begin{equation}\label{6}
\rho (\mathbb{D}\mathbb{Q}(G)) \geq {\frac{{3\mathbb{D}_n - 1+ \sqrt
{(\mathbb{D}_n+1)^2 + 8\sum\nolimits_{k = 1}^{n - 1} {(\mathbb{D}_k
- \mathbb{D}_n )} } }}{2}} \cite{Duan2013}.
\end{equation}

In this paper, using a positive scale vector and their
intersections, we give some upper and lower bounds for the spectral
radius of a nonnegative irreducible matrix and characterize the
equality cases. These bounds theoretically improve and generalize
some known results in \cite{Duan2013}. Applying these bounds to
various matrices associated with a graph, including the adjacency
matrix, the signless Laplacian matrix, the distance matrix, the
distance signless Laplacian matrix, we obtain some new upper and
lower bounds on spectral radius of various matrices of a graph as
mentioned above, which generalize and improve some known results.

\section*{2. Upper and lower bounds for spectral radius of a nonnegative matrix}

In \cite{Tian2009}, Tian et al. presented some new inclusion
intervals of matrix singular values. These intervals are based
mainly on the use of positive scale vectors and their intersections.
In this section, using a similar technique, we shall give some upper
and
lower bounds for spectral radius of a nonnegative matrix.\\
\\
\textbf{Theorem 1.} {\it Let $A=(a_{ij})$ be an $n\times n$
nonnegative irreducible matrix and $c=(c_1,c_2,\ldots,c_n)^T$ be any
vector with positive components. For $1\leq i\leq n$, take
\[
M_i = \frac{1}{{c_i }}\sum\limits_{j = 1}^n {a_{ij}c_j } ,\; M =
\mathop {\max }\limits_i \{ a_{ii} \} ,\; N = \mathop {\max
}\limits_{i \ne j} \{ \frac{{a_{ij} c_j }}{{c_i }}\}.
\]
Assume that $M_1\geq M_2\geq\cdots\geq M_n$ and $N>0$. Then, for
$1\leq i\leq n$,
\begin{equation}\label{7}
\rho (A) \le {\frac{{M_i +M-N  + \sqrt {(M_i  - M  + N)^2 +
4N\sum\nolimits_{k = 1}^{i - 1} {(M_k  - M_i )} } }}{2}}.
\end{equation}
Moreover, the equality holds in (\ref{7}) if and only if $M_1=
M_2=\cdots= M_n$ or for some $2\leq t\leq i$, $A$ satisfies the
following conditions:

(i) $a_{kk}=M$ for $1\leq k\leq t-1$,

(ii) $\frac{a_{kl}c_l}{c_k}=N$ for $1\leq k\leq n$, $1\leq l\leq
t-1$ with $k\neq l$,

(iii) $M_t=\cdots= M_n$.
}\\
\\
\textbf{Proof.} Let $U=\text{diag}(c_1,c_2,\ldots, c_n)$. Obviously,
$A$ and $B=U^{-1}AU$ have the same eigenvalues. Thus, applying
Theorem 2.1 in \cite{Duan2013} to $B$, one may obtain the required result. $\Box$\\
\\
\textbf{Corollary 1}\cite{Duan2013}. {\it Let $A=(a_{ij})$ be an
$n\times n$ nonnegative irreducible matrix with row sums $r_1\geq
r_2\geq\cdots\geq r_n$. Also let $ M = \mathop {\max
}\nolimits_{1\leq i\leq n}\{ a_{ii} \} ,\; N = \mathop {\max
}\nolimits_{i \ne j} \{ a_{ij}\}$. Assume that $N>0$. Then, for
$1\leq i\leq n$,
\begin{equation}\label{8}
\rho (A) \le {\frac{{r_i +M-N  + \sqrt {(r_i  - M  + N)^2 +
4N\sum\nolimits_{k = 1}^{i - 1} {(r_k  - r_i )} } }}{2}}.
\end{equation}
Moreover, the equality holds in (\ref{8}) if and only if $r_1=
r_2=\cdots= r_n$ or for some $2\leq t\leq i$, $A$ satisfies the
following conditions:

(i) $a_{kk}=M$ for $1\leq k\leq t-1$,

(ii) $a_{kl}=N$ for $1\leq k\leq n$, $1\leq l\leq t-1$ with $k\neq
l$,

(iii) $r_t=\cdots= r_n$.}\\
\\
\textbf{Proof.} Take $c=(1,1,\ldots,1)^T$ in Theorem 1, one may
obtain the required result. $\Box$\\
\\
\textbf{Corollary 2}\cite{Chen2013}. {\it Let $A=(a_{ij})$ be an
$n\times n$ nonnegative irreducible matrix with diagonal entries 0,
its row sums $r_1\geq r_2\geq\cdots\geq r_n$. Also let $ N = \mathop
{\max }\nolimits_{i \ne j} \{ a_{ij}\}$. Assume that $N>0$. Then,
for $1\leq i\leq n$,
\begin{equation}\label{9}
\rho (A) \le {\frac{{r_i-N  + \sqrt {(r_i+ N)^2 + 4N\sum\nolimits_{k
= 1}^{i - 1} {(r_k  - r_i )} } }}{2}}.
\end{equation}
Moreover, the equality holds in (\ref{9}) if and only if $r_1=
r_2=\cdots= r_n$ or for some $2\leq t\leq i$, $A$ satisfies the
following conditions:

(i) $a_{kl}=N$ for $1\leq k\leq n$, $1\leq l\leq t-1$ with $k\neq
l$,

(ii) $r_t=\cdots= r_n$.}\\
\\
\textbf{Proof.} Take $c=(1,1,\ldots,1)^T$ and $M=0$ in Theorem 1,
one may obtain the required result. $\Box$\\
\\
\textbf{Remark 1.} From Corollaries 1 and 2, one may see that
Theorem 1 improves and generalizes some results in
\cite{Chen2013,Duan2013}. On the other hand, we always assume that
$M_1\geq M_2\geq\cdots\geq M_n$ in Theorem 1 because the matrix
$B=U^{-1}AU$ in the proof of Theorem 1 may satisfy this condition
under permutation similarity. Finally, Chen et al.\cite{Chen2013}
claim that the equality holds in (\ref{9}) if and only if $r_1=
r_2=\cdots= r_n$ or for some $2\leq i\leq n$, $A$ satisfies
$r_1=r_2=\cdots=r_{i-1}>r_i=\cdots= r_n$ and $a_{kl}=N$ for $1\leq
k\leq n$, $1\leq l\leq i-1$ with $k\neq l$ (see Theorem 1.9 in
\cite{Chen2013}). However, it seems to be inaccurate. For example,
let
\[
A = \left( {\begin{array}{*{20}c}
   0 & 4 & 2 & 3 & 3  \\
   4 & 0 & 2 & 2 & 3  \\
   4 & 4 & 0 & 1 & 1  \\
   4 & 4 & 1 & 0 & 1  \\
   4 & 4 & 1 & 1 & 0  \\
\end{array}} \right).
\]
For $i=3$, applying the inequality (\ref{9}), one gets $\rho (A) \le
\frac{{6 + \sqrt {244} }}{2} \approx 10.8102$. In fact, by direct
calculation, one has $\rho (A) \approx 10.8102$. Hence the equality
in (\ref{9}) holds for the matrix $A$. But
$r_1=12,r_2=11,r_3=r_4=r_5=10$.\\
\\
\textbf{Theorem 2.} {\it Let $A=(a_{ij})$ be an $n\times n$
nonnegative irreducible matrix and $c=(c_1,c_2,\ldots,c_n)^T$ be any
vector with positive components. For $1\leq i\leq n$, take
\[
M_i = \frac{1}{{c_i }}\sum\limits_{j = 1}^n {a_{ij}c_j } ,\; S =
\mathop {\min }\limits_i \{ a_{ii} \} ,\; T = \mathop {\min
}\limits_{i \ne j} \{ \frac{{a_{ij} c_j }}{{c_i }}\}.
\]
Assume that $M_1\geq M_2\geq\cdots\geq M_n$. Then
\begin{equation}\label{10}
\rho (A) \geq {\frac{{M_n +S-T + \sqrt {(M_n  - S  + T)^2 +
4T\sum\nolimits_{k = 1}^{n- 1} {(M_k  - M_n )} } }}{2}}.
\end{equation}
Moreover, the equality holds in (\ref{10}) if and only if $M_1=
M_2=\cdots= M_n$ or $T>0$, for some $2\leq t\leq n$, $A$ satisfies
the following conditions:

(i) $a_{kk}=S$ for $1\leq k\leq t-1$,

(ii) $\frac{a_{kl}c_l}{c_k}=T$ for $1\leq k\leq n$, $1\leq l\leq
t-1$ with $k\neq l$,

(iii) $M_t=\cdots= M_n$.
}\\
\\
\textbf{Proof.} Let $U=\text{diag}(c_1,c_2,\ldots,c_n)$. Obviously,
$A$ and $B=U^{-1}AU$ have the same eigenvalues. Thus, applying
Theorem 2.2 in \cite{Duan2013} to $B$, one may obtain the required result. $\Box$\\
\\
\textbf{Corollary 3}\cite{Duan2013}. {\it Let $A=(a_{ij})$ be an
$n\times n$ nonnegative irreducible matrix with row sums $r_1\geq
r_2\geq\cdots\geq r_n$. Also let $ M = \mathop {\min
}\nolimits_{1\leq i\leq n}\{ a_{ii} \} ,\; N = \mathop {\min
}\nolimits_{i \ne j} \{ a_{ij}\}$. Then
\begin{equation}\label{11}
\rho (A) \geq {\frac{{r_n +S-T + \sqrt {(r_n  - S  + T)^2 +
4T\sum\nolimits_{k = 1}^{n- 1} {(r_k  - r_n )} } }}{2}}.
\end{equation}
Moreover, the equality holds in (\ref{11}) if and only if $r_1=
r_2=\cdots= r_n$ or $T>0$, for some $2\leq t\leq n$, $A$ satisfies
the following conditions:

(i) $a_{kk}=S$ for $1\leq k\leq t-1$,

(ii) $a_{kl}=T$ for $1\leq k\leq n$, $1\leq l\leq t-1$ with $k\neq
l$,

(iii) $r_t=\cdots= r_n$.}\\
\\
\textbf{Proof.} Take $c=(1,1,\ldots,1)^T$ in Theorem 1, one may
obtain the required result. $\Box$\\
\\
\textbf{Remark 2.} Corollary 3 implies that Theorem 2 is a
generalization and improvement on Theorem 2.2 in \cite{Duan2013}.

\section*{3. Bounds for spectral radius of various matrices associated with a graph}

Let $G=(V,E)$ be a simple graph with vertex set
$V=\{v_1,v_2,\ldots,v_n\}$ and edge set $E$. Let $d_1\geq d_2\geq
\ldots \geq d_n$ be degree sequence of $G$. Denote by $m_i$ the
average degree of a vertex $v_i$, which is $\frac{\sum_{j\sim
i}d_j}{d_i}$. In \cite{Liu2010, Tian2012}, authors introduced the
following notation: for any real number $\alpha$,
\[
({}^\alpha m)_i  = \frac{{\sum\nolimits_{j \sim i} {d_j^\alpha  }
}}{{d_i^\alpha  }},
\]
which is called the generalized average degree of $v_i$. Note that
$d_i=(^0m)_i$ and $m_i=(^1m)_i$. Moreover, they also presented some
upper and lower bounds for the (Laplacian) spectral radius of $G$
with parameter $\alpha$. Next, using a similar method, we also give
some bounds for spectral radius of various matrices associated
with $G$.\\
\\
\textbf{Theorem 3.} {\it Let $G$ be a connected graph on $n\geq 2$
vertices with generalized average degree $(^\alpha m)_1\geq (^\alpha
m)_2\geq\cdots\geq(^\alpha m)_n$. Then, for $1\leq i\leq n$,
\begin{equation}\label{12}
\rho (A(G)) \le {\frac{{(^\alpha m)_i -N  + \sqrt {((^\alpha m)_i +
N)^2 + 4N\sum\nolimits_{k = 1}^{i - 1} {((^\alpha m)_k  - (^\alpha
m)_i )} } }}{2}},
\end{equation}
where $N = \mathop {\max }\nolimits_{i \sim j} \{ {d_j^\alpha
\mathord{\left/
 {\vphantom {d_j^\alpha d_i^\alpha}} \right.
 \kern-\nulldelimiterspace} d_i^\alpha}\}$. Moreover, the equality holds in (\ref{12}) if and only if
$(^\alpha m)_1= (^\alpha m)_2=\cdots=(^\alpha m)_n$ or for some
$2\leq t\leq i$, $G$ has one of the following properties:

(i) If $\alpha=0$, then $G$ is a bidegreed graph in which $d_1= d_2=
\cdots =d_{t-1}=n-1>d_t= d_{t+1}= \cdots =d_n$;

(ii) If $\alpha> 0$, then $G$ is a bidegreed graph in which
$(^\alpha m)_1> (^\alpha
m)_2=\cdots=(^\alpha m)_n$ and $d_1=n-1>d_2= \cdots =d_n$.}\\
\\
\textbf{Proof.} Let $c=(d_1^\alpha,d_2^\alpha,\ldots,d_n^\alpha)^T$
in (\ref{7}). Obviously, $c$ is a positive vector as $G$ is
connected. Now apply Theorem 1 to $A(G)$. Notice that $M_i=(^\alpha
m)_i, M=0, N = \mathop {\max }\nolimits_{i \sim j} \{ {d_j^\alpha
\mathord{\left/
 {\vphantom {d_j^\alpha d_i^\alpha}} \right.
 \kern-\nulldelimiterspace} d_i^\alpha}\}$. We readily get the required
 upper bound (\ref{12}). Again from Theorem 1, the the equality holds in (\ref{12}) if and only if
$(^\alpha m)_1= (^\alpha m)_2=\cdots=(^\alpha m)_n$ or for some
$2\leq t\leq i$, $A(G)=(a_{ij})$ satisfies the following conditions:

(a) $\frac{a_{kl}d_l^\alpha}{d_k^\alpha}=N$ for $1\leq k\leq n$,
$1\leq l\leq t-1$ with $k\neq l$,

(b) $(^\alpha m)_t=\cdots=(^\alpha m)_n$.

Consider the following three cases: $\alpha=0$, $\alpha>0$ and
$\alpha<0$. If $\alpha=0$, then the above (a) implies that $d_1=
d_2= \cdots =d_{t-1}=n-1$. Since $(^0 m)_i=d_i$. the above (b)
becomes $d_t= d_{t+1}= \cdots =d_n$. If $\alpha>0$, then the above
(a) implies $d_1=n-1$ and $d_2=d_3=\cdots=d_n$ as $N\geq 1$ is a
constant. Thus $(^\alpha m)_1\geq (^\alpha m)_2=\cdots=(^\alpha
m)_n$. Finally, consider $\alpha<0$. From the (a) and $N\geq 1$, one
gets $d_1=n-1\leq d_2=\cdots=d_n$, which implies that $G$ is a
complete graph.

This completes our proof. $\Box$\\
\\
\textbf{Remark 3.} Take $\alpha=0$ in (\ref{12}), we obtain the
upper bound (\ref{1}). Hence the upper bound (12) improves and
generalizes some results in \cite{Duan2013,Liu2013}.\\
\\
\textbf{Corollary 4}\cite{Huang20xx}. {\it Let $G$ be a connected
graph on $n\geq 2$ vertices with average degree $ m_1\geq
m_2\geq\cdots\geq m_n$. Then, for $1\leq i\leq n$,
\begin{equation}\label{13}
\rho (A(G)) \le {\frac{{ m_i -N  + \sqrt {(m_i + N)^2 +
4N\sum\nolimits_{k = 1}^{i - 1} {(m_k  - m_i )} } }}{2}},
\end{equation}
where $N = \mathop {\max }\nolimits_{i \sim j} \{ {d_j
\mathord{\left/
 {\vphantom {d_j d_i}} \right.
 \kern-\nulldelimiterspace} d_i}\}$. Moreover, the equality holds in (\ref{13}) if and only if
$ m_1= m_2=\cdots= m_n$, i.e., $G$ is pseudo-regular.}\\
\\
\textbf{Proof.} Take $\alpha=1$ in (\ref{12}), we obtain the upper
bound (\ref{13}) as $(^1m)_i=m_i$. Moreover, the equality holds in
(\ref{13}) if and only if $m_1= m_2=\cdots=m_n$ or $G$ is a
bidegreed graph in which $m_1>m_2=\cdots=m_n$ and $d_1=n-1>d_2=
\cdots =d_n$. If $G$ is a bidegreed graph in which
$m_1>m_2=\cdots=m_n$ and $d_1=n-1>d_2= \cdots =d_n=\delta$, then
$m_1=\delta$ and $m_2=\delta-1+\frac{n-1}{\delta}$. $m_1>m_2$
implies $\delta>n-1$, which is impossible. Hence, $m_1=
m_2=\cdots=m_n$. $\Box$\\
\\
\textbf{Theorem 4.} {\it Let $G$ be a connected graph on $n\geq 2$
vertices with $(^\alpha m)_1+d_1\geq (^\alpha
m)_2+d_2\geq\cdots\geq(^\alpha m)_n+d_n$. Then, for $1\leq i\leq n$,
\begin{equation}\label{14}
\rho (Q(G)) \le {\frac{{(^\alpha m)_i+d_i+\Delta -N  + \sqrt
{((^\alpha m)_i+d_i-\Delta + N)^2 + 4N\sum\nolimits_{k = 1}^{i - 1}
{((^\alpha m)_k +d_k- (^\alpha m)_i -d_i)} } }}{2}},
\end{equation}
where $\Delta$ is the maximum degree of $G$, $N = \mathop {\max
}\nolimits_{i \sim j} \{ {d_j^\alpha \mathord{\left/
 {\vphantom {d_j^\alpha d_i^\alpha}} \right.
 \kern-\nulldelimiterspace} d_i^\alpha}\}$. Moreover, the equality holds in (\ref{14}) if and only if
$(^\alpha m)_1+d_1= (^\alpha m)_2+d_2=\cdots=(^\alpha m)_n+d_n$ or
for some $2\leq t\leq i$, $G$ has one of the following properties:

(i) If $\alpha=0$, then $G$ is a bidegreed graph in which $d_1= d_2=
\cdots =d_{t-1}=n-1>d_t= d_{t+1}= \cdots =d_n$;

(ii) If $\alpha> 0$, then $G$ is a bidegreed graph in which
$(^\alpha m)_1+d_1> (^\alpha
m)_2+d_2=\cdots=(^\alpha m)_n+d_n$ and $d_1=n-1>d_2= \cdots =d_n$.}\\
\\
\textbf{Proof.} Let $c=(d_1^\alpha,d_2^\alpha,\ldots,d_n^\alpha)^T$.
Obviously, $c$ is a positive vector as $G$ is connected. Now apply
Theorem 1 to $Q(G)$. Notice that $M_i=(^\alpha m)_i+d_i, M=\Delta, N
= \mathop {\max }\nolimits_{i \sim j} \{ {d_j^\alpha \mathord{\left/
 {\vphantom {d_j^\alpha d_i^\alpha}} \right.
 \kern-\nulldelimiterspace} d_i^\alpha}\}$. Thus we obtain the required
 upper bound (\ref{14}). Again from Theorem 1, the the equality in (\ref{14}) holds if and only if
$(^\alpha m)_1+d_1= (^\alpha m)_2+d_2=\cdots=(^\alpha m)_n+d_n$ or
for some $2\leq t\leq i$, $Q(G)=(q_{ij})$ satisfies the following
conditions:

(a) $q_{kk}=\Delta$ for $1\leq k\leq t-1$,

(b) $\frac{q_{kl}d_l^\alpha}{d_k^\alpha}=N$ for $1\leq k\leq n$,
$1\leq l\leq t-1$ with $k\neq l$,

(c) $(^\alpha m)_t+d_t=\cdots=(^\alpha m)_n+d_n$.

Consider the following three cases: $\alpha=0$, $\alpha>0$ and
$\alpha<0$. If $\alpha=0$, then the above (b) implies $d_1= d_2=
\cdots =d_{t-1}=n-1$. Since $(^0 m)_i=d_i$. the above (c) becomes
$d_t= d_{t+1}= \cdots =d_n$. If $\alpha>0$, then the above (b) shows
that $d_1=n-1$ and $d_2=d_3=\cdots=d_n$ as $N\geq 1$ is a constant.
Thus $(^\alpha m)_1+d_1\geq (^\alpha m)_2+d_2=\cdots=(^\alpha
m)_n+d_n$. Finally, consider $\alpha<0$. From the (b) and $N\geq 1$,
one has $d_1=n-1\leq d_2=\cdots=d_n$, which implies that $G$ is a
complete graph.

This completes our proof. $\Box$\\
\\
\textbf{Remark 4.} Take $\alpha=0$ in (\ref{14}), we may obtain the
upper bound (\ref{2}). Hence the upper bound (\ref{14}) improves and
generalizes some results in \cite{Duan2013,Cui2013}.\\
\
\\
\textbf{Corollary 5.} {\it Let $G$ be a connected graph on $n\geq 2$
vertices with $m_1+d_1\geq m_2+d_2\geq\cdots\geq m_n+d_n$. Then, for
$1\leq i\leq n$,
\begin{equation}\label{15}
\rho (Q(G)) \le {\frac{{m_i+d_i+\Delta -N  + \sqrt {(m_i+d_i-\Delta
+ N)^2 + 4N\sum\nolimits_{k = 1}^{i - 1} {(m_k +d_k- m_i -d_i)} }
}}{2}},
\end{equation}
where $\Delta$ is the maximum degree of $G$, $N = \mathop {\max
}\nolimits_{i \sim j} \{ {d_j \mathord{\left/
 {\vphantom {d_j d_i}} \right.
 \kern-\nulldelimiterspace} d_i}\}$. Moreover, the equality holds in (\ref{15}) if and only if
$m_1+d_1= m_2+d_2=\cdots=m_n+d_n$ or for some $2\leq i\leq n$, $G$
is a bidegreed graph in
which $m_1+d_1> m_2+d_2=\cdots=m_n+d_n$ and $d_1=n-1>d_2= \cdots =d_n$.}\\
\\
\textbf{Proof.} Take $\alpha=1$ in (\ref{14}), the required result
follows as $(^1m)_i=m_i$. $\Box$\\

Remark that the equality holds in (\ref{15}) for some $2\leq i\leq
n$ whenever $G$ is a bidegreed graph in which $m_1+d_1>
m_2+d_2=\cdots=m_n+d_n$ and $d_1=n-1>d_2= \cdots =d_n$. For example,
consider the following graph $G$ shown in Fig.1. For $i=2$, applying
the inequality (\ref{15}), one gets $\rho (Q(G)) \le \frac{{7 +
\sqrt {17} }}{2} \approx 5.5616$. By direct calculation, one has
$\rho (Q(G)) \approx 5.5616$. Hence the equality in (\ref{15}) holds
for the graph $G$ shown in Fig.1. In fact, $d_1=4>d_2=\cdots=d_5=2$
and $m_1+d_1=6>m_2+d_2=\cdots=m_5+d_5=5$.
\begin{figure} [h t
b p] \centering
%\begin{minipage}[c]{12cm}
%\centering
\includegraphics[width=5cm]{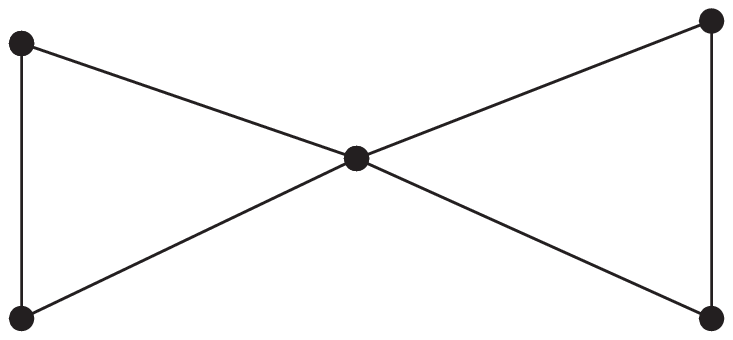}\\
\centerline{}
%\end{minipage}\\
\centerline{\small Fig.1. The shown graph $G$}
\end{figure}\\

Let $G=(V,E)$ be a simple graph with vertex set
$V=\{v_1,v_2,\ldots,v_n\}$ and edge set $E$. For any vertex $v_i\in
V$, denote the transmission of $v_i$ by $\mathbb{D}_i$. Denote by
$\mathbb{M}_i$ the \emph{average transmission} of a vertex $v_i$,
which is $\frac{\sum_{j=1}^nd_{ij}\mathbb{D}_j}{\mathbb{D}_i}$. Now
we introduce the following notation: for any real number $\alpha$,
\[
({}^\alpha \mathbb{M})_i  =
\frac{\sum_{j=1}^nd_{ij}\mathbb{D}_j^\alpha}{\mathbb{D}_i^\alpha},
\]
which is called the \emph{generalized average transmission} of
$v_i$. Note that
$\mathbb{D}_i=(^0\mathbb{M})_i$ and $\mathbb{M}_i=(^1\mathbb{M})_i$.\\
\\
\textbf{Theorem 5.} {\it Let $G$ be a connected graph on $n\geq2$
vertices with generalized average transmissions
$(^\alpha\mathbb{M})_1\geq(^\alpha\mathbb{M})_2\geq\cdots\geq(^\alpha\mathbb{M})_n$.
Then, for $1\leq i\leq n$,
\begin{equation}\label{16}
\rho (\mathbb{D}(G)) \le {\frac{{(^\alpha \mathbb{M})_i -N  + \sqrt
{((^\alpha \mathbb{M})_i + N)^2 + 4N\sum\nolimits_{k = 1}^{i - 1}
{((^\alpha \mathbb{M})_k  - (^\alpha \mathbb{M})_i )} } }}{2}},
\end{equation}
where $N = \mathop {\max }\nolimits_{1\leq i,j\leq n} \{
{d_{ij}\mathbb{D}_j^\alpha \mathord{\left/
 {\vphantom {d_{ij}\mathbb{D}_j^\alpha \mathbb{D}_i^\alpha}} \right.
 \kern-\nulldelimiterspace} \mathbb{D}_i^\alpha}\}$. Moreover, the equality holds in (\ref{16}) if and only if
$(^\alpha \mathbb{M})_1= (^\alpha \mathbb{M})_2=\cdots=(^\alpha
\mathbb{M})_n$ or for some $2\leq t\leq i$, $\mathbb{D}(G)=(d_{ij})$
satisfies the following conditions:

(i) $\frac{d_{kl}\mathbb{D}_l^\alpha}{\mathbb{D}_k^\alpha}=N$ for
$1\leq k\leq n$, $1\leq l\leq t-1$ with $k\neq l$,

(ii) $(^\alpha \mathbb{M})_t= (^\alpha
\mathbb{M})_{t+1}=\cdots=(^\alpha \mathbb{M})_n$.
}\\
\\
\textbf{Proof.} Let
$c=(\mathbb{D}_1^\alpha,\mathbb{D}_2^\alpha,\ldots,\mathbb{D}_n^\alpha)^T$.
Apply Theorem 1 to $\mathbb{D}(G)$. Notice that $M_i=(^\alpha
\mathbb{M})_i, M=0, N = \mathop {\max }\nolimits_{1\leq i,j\leq n}
\{ {d_{ij}\mathbb{D}_j^\alpha \mathord{\left/
 {\vphantom {d_{ij}\mathbb{D}_j^\alpha \mathbb{D}_i^\alpha}} \right.
 \kern-\nulldelimiterspace} \mathbb{D}_i^\alpha}\}$. We easily obtain the required
result.\\
\\
\textbf{Corollary 6}\cite{Duan2013}. {\it Let $G$ be a connected
graph on $n\geq 2$ vertices with transmissions $\mathbb{D}_1\geq
\mathbb{D}_2\geq \ldots \geq \mathbb{D}_n$. Also let $d$ be the
diameter of $G$. Then, for $1\leq i\leq n$, the inequality (\ref{3})
holds. Moreover, the equality holds in (\ref{3}) if and only if $
\mathbb{D}_1=  \mathbb{D}_2=\cdots= \mathbb{D}_n$.
}\\
\\
\textbf{Proof.} Take $\alpha=0$ in (\ref{16}), the upper bound
(\ref{3}) follows. Again from Theorem 5, the equality holds in
(\ref{3}) if and only if $ \mathbb{D}_1=  \mathbb{D}_2=\cdots=
\mathbb{D}_n$ or for some $2\leq t\leq i$, $\mathbb{D}(G)=(d_{ij})$
satisfies the following conditions:

(a) $d_{kl}=N=d$ for $1\leq k\leq n$, $1\leq l\leq t-1$ with $k\neq
l$,

(b) $\mathbb{D}_t= \mathbb{D}_{t+1}=\cdots=\mathbb{D}_n$\\
It follows from above (a) that $d=1$. Hence $G$ is a complete graph
for the latter case. $\Box$\\

If we take $\alpha=1$ in (\ref{16}), then the following corollary is
obtained:\\
\\
\textbf{Corollary 7.} {\it Let $G$ be a connected graph on $n\geq 2$
vertices with average transmissions $ \mathbb{M}_1\geq
\mathbb{M}_2\geq\cdots\geq \mathbb{M}_n$. Then, for $1\leq i\leq n$,
\begin{equation}\label{17}
\rho (\mathbb{D}(G)) \le {\frac{{ \mathbb{M}_i -N  + \sqrt
{(\mathbb{M}_i + N)^2 + 4N\sum\nolimits_{k = 1}^{i - 1}
{(\mathbb{M}_k  - \mathbb{M}_i )} } }}{2}},
\end{equation}
where $N = \mathop {\max }\nolimits_{1\leq i,j\leq n} \{
{d_{ij}\mathbb{D}_j \mathord{\left/
 {\vphantom {d_{ij}\mathbb{D}_j \mathbb{D}_i}} \right.
 \kern-\nulldelimiterspace} \mathbb{D}_i}\}$. Moreover, the equality holds in (\ref{17}) if and only if
$ \mathbb{M}_1=  \mathbb{M}_2=\cdots= \mathbb{M}_n$ or for some
$2\leq t\leq i$, $\mathbb{D}(G)=(d_{ij})$ satisfies the following
conditions:

(i) $\frac{d_{kl}\mathbb{D}_l}{\mathbb{D}_k}=N$ for $1\leq k\leq n$,
$1\leq l\leq t-1$ with $k\neq l$,

(ii) $\mathbb{M}_t=
\mathbb{M}_{t+1}=\cdots=\mathbb{M}_n$. }\\
\\
\textbf{Theorem 6.} {\it Let $G$ be a connected graph on $n\geq2$
vertices with generalized average transmissions
$(^\alpha\mathbb{M})_1\geq(^\alpha\mathbb{M})_2\geq\cdots\geq(^\alpha\mathbb{M})_n$.
Then,
\begin{equation}\label{18}
\rho (\mathbb{D}(G)) \geq {\frac{{(^\alpha \mathbb{M})_n -T  + \sqrt
{((^\alpha \mathbb{M})_n + T)^2 + 4T\sum\nolimits_{k = 1}^{n - 1}
{((^\alpha \mathbb{M})_k  - (^\alpha \mathbb{M})_n )} } }}{2}},
\end{equation}
where $T = \mathop {\min }\nolimits_{1\leq i\neq j\leq n} \{
{d_{ij}\mathbb{D}_j^\alpha \mathord{\left/
 {\vphantom {d_{ij}\mathbb{D}_j^\alpha \mathbb{D}_i^\alpha}} \right.
 \kern-\nulldelimiterspace} \mathbb{D}_i^\alpha}\}$. Moreover, the equality holds in (\ref{18}) if and only if
$(^\alpha \mathbb{M})_1= (^\alpha \mathbb{M})_2=\cdots=(^\alpha
\mathbb{M})_n$ or for some $2\leq t\leq n$, $\mathbb{D}(G)=(d_{ij})$
satisfies the following conditions:

(i) $\frac{d_{kl}\mathbb{D}_l^\alpha}{\mathbb{D}_k^\alpha}=T$ for
$1\leq k\leq n$, $1\leq l\leq t-1$ with $k\neq l$,

(ii) $(^\alpha \mathbb{M})_t= (^\alpha
\mathbb{M})_{t+1}=\cdots=(^\alpha \mathbb{M})_n$. }\\
\\
\textbf{Proof.} Apply Theorem 2 to $\mathbb{D}(G)$. The rest proof
is similar to that Theorem 5, omitted. $\Box$\\
\\
\textbf{Remark 5.} Take $\alpha=0$ in (\ref{18}), we may obtain the
bound (\ref{5}). Hence the lower bound (\ref{18}) improves and
generalizes the result in \cite{Duan2013}. If $\alpha=1$ in
(\ref{18}), then
\begin{equation*}
\rho (\mathbb{D}(G)) \geq {\frac{{ \mathbb{M}_n -T  + \sqrt
{(\mathbb{M}_n + T)^2 + 4T\sum\nolimits_{k = 1}^{n - 1} {(
\mathbb{M}_k  -  \mathbb{M}_n) } } }}{2}},
\end{equation*}
where $T = \mathop {\min }\nolimits_{1\leq i\neq j\leq n} \{
{d_{ij}\mathbb{D}_j \mathord{\left/
 {\vphantom {d_{ij}\mathbb{D}_j \mathbb{D}_i}} \right.
 \kern-\nulldelimiterspace} \mathbb{D}_i}\}$.
\\

Apply Theorem 1 to the distance signless Laplacian matrix
$\mathbb{D}\mathbb{Q}(G)$, we easily arrive at:\\
\\
\textbf{Theorem 7.} {\it Let $G$ be a connected graph on $n\geq 2$
vertices with $(^\alpha \mathbb{M})_1+\mathbb{D}_1\geq
(^\alpha\mathbb{M} )_2+\mathbb{D}_2\geq\cdots\geq(^\alpha
\mathbb{M})_n+\mathbb{D}_n$. Then, for $1\leq i\leq n$,
\begin{equation}\label{19}
\small \rho (\mathbb{D}\mathbb{Q}(G)) \le {\frac{{(^\alpha
\mathbb{M})_i+\mathbb{D}_i+M -N  + \sqrt {((^\alpha
\mathbb{M})_i+\mathbb{D}_i-M + N)^2 + 4N\sum\nolimits_{k = 1}^{i -
1} {((^\alpha \mathbb{M})_k +\mathbb{D}_k- (^\alpha \mathbb{M})_i
-\mathbb{D}_i)} } }}{2}},
\end{equation}
where $M=\max\nolimits_{1\leq i\leq n}\{\mathbb{D}_i\}$, $N =
\mathop {\max }\nolimits_{1\leq i\neq j\leq n} \{
{d_{ij}\mathbb{D}_j^\alpha \mathord{\left/
 {\vphantom {d_{ij}\mathbb{D}_j^\alpha \mathbb{D}_i^\alpha}} \right.
 \kern-\nulldelimiterspace} \mathbb{D}_i^\alpha}\}$. Moreover, the equality holds in (\ref{19}) if and only if
$(^\alpha \mathbb{M})_1+\mathbb{D}_1= (^\alpha
\mathbb{M})_2+\mathbb{D}_2=\cdots=(^\alpha
\mathbb{M})_n+\mathbb{D}_n$ or for some $2\leq t\leq i$,
$\mathbb{D}\mathbb{Q}(G)=(d_{ij})$ satisfies the following
conditions:

(i) $d_{kk}=M$ for $1\leq k\leq t-1$,

(ii) $\frac{d_{kl}\mathbb{D}_l^\alpha}{\mathbb{D}_k^\alpha}=N$ for
$1\leq k\leq n$, $1\leq l\leq t-1$ with $k\neq l$,

(iii) $(^\alpha \mathbb{M})_t+\mathbb{D}_t= (^\alpha
\mathbb{M})_{t+1}+\mathbb{D}_{t+1}=\cdots=(^\alpha \mathbb{M})_n+\mathbb{D}_n$. }\\

If we apply Theorem 2 to the distance signless Laplacian matrix
$\mathbb{D}\mathbb{Q}(G)$, then the following theorem is obtained.\\
\\
\textbf{Theorem 8.} {\it Let $G$ be a connected graph on $n\geq 2$
vertices with $(^\alpha \mathbb{M})_1+\mathbb{D}_1\geq
(^\alpha\mathbb{M} )_2+\mathbb{D}_2\geq\cdots\geq(^\alpha
\mathbb{M})_n+\mathbb{D}_n$. Then
\begin{equation}\label{20}
\small \rho (\mathbb{D}\mathbb{Q}(G)) \geq {\frac{{(^\alpha
\mathbb{M})_n+\mathbb{D}_n+S -T  + \sqrt {((^\alpha
\mathbb{M})_n+\mathbb{D}_n-S + T)^2 + 4T\sum\nolimits_{k = 1}^{n -
1} {((^\alpha \mathbb{M})_k +\mathbb{D}_k- (^\alpha \mathbb{M})_n
-\mathbb{D}_n)} } }}{2}},
\end{equation}
where $S=\min\nolimits_{1\leq i\leq n}\{\mathbb{D}_i\}$, $T= \mathop
{\min}\nolimits_{1\leq i\neq j\leq n} \{ {d_{ij}\mathbb{D}_j^\alpha
\mathord{\left/
 {\vphantom {d_{ij}\mathbb{D}_j^\alpha \mathbb{D}_i^\alpha}} \right.
 \kern-\nulldelimiterspace} \mathbb{D}_i^\alpha}\}$. Moreover, the
equality holds in (\ref{20}) if and only if $(^\alpha
\mathbb{M})_1+\mathbb{D}_1= (^\alpha
\mathbb{M})_2+\mathbb{D}_2=\cdots=(^\alpha
\mathbb{M})_n+\mathbb{D}_n$ or for some $2\leq t\leq n$,
$\mathbb{D}\mathbb{Q}(G)=(d_{ij})$ satisfies the following
conditions:

(i) $d_{kk}=S$ for $1\leq k\leq t-1$,

(ii) $\frac{d_{kl}\mathbb{D}_l^\alpha}{\mathbb{D}_k^\alpha}=T$ for
$1\leq k\leq n$, $1\leq l\leq t-1$ with $k\neq l$,

(iii) $(^\alpha \mathbb{M})_t+\mathbb{D}_t= (^\alpha
\mathbb{M})_{t+1}+\mathbb{D}_{t+1}=\cdots=(^\alpha
\mathbb{M})_n+\mathbb{D}_n$.\\
  }\
\\
\textbf{Remark 6.} If $\alpha=0$ in (\ref{19}), then the bound upper
(\ref{4}) is obtained. Hence the upper bound (\ref{19}) improves and
generalizes the result in \cite{Duan2013}. If $\alpha=1$ in
(\ref{19}), then
\begin{equation*}
\small \rho (\mathbb{D}\mathbb{Q}(G)) \le {\frac{{
\mathbb{M}_i+\mathbb{D}_i+M -N  + \sqrt {(
\mathbb{M}_i+\mathbb{D}_i-M + N)^2 + 4N\sum\nolimits_{k = 1}^{i - 1}
{( \mathbb{M}_k +\mathbb{D}_k- \mathbb{M}_i -\mathbb{D}_i)} }
}}{2}},
\end{equation*}
where $M=\max\nolimits_{1\leq i\leq n}\{\mathbb{D}_i\}$, $N =
\mathop {\max }\nolimits_{1\leq i\neq j\leq n} \{
{d_{ij}\mathbb{D}_j^\alpha \mathord{\left/
 {\vphantom {d_{ij}\mathbb{D}_j^\alpha \mathbb{D}_i^\alpha}} \right.
 \kern-\nulldelimiterspace} \mathbb{D}_i^\alpha}\}$. Similarly, if
$\alpha=0$ in (\ref{20}), then the lower upper (\ref{6}) follows.
Hence the lower bound (\ref{20}) improves and generalizes the result
in \cite{Duan2013}. If $\alpha=1$ in (\ref{20}), then
\begin{equation*}
\small \rho (\mathbb{D}\mathbb{Q}(G)) \geq {\frac{{
\mathbb{M}_n+\mathbb{D}_n+S -T  + \sqrt {(
\mathbb{M}_n+\mathbb{D}_n-S + T)^2 + 4T\sum\nolimits_{k = 1}^{n - 1}
{( \mathbb{M}_k +\mathbb{D}_k- \mathbb{M}_n -\mathbb{D}_n)} }
}}{2}},
\end{equation*}
where $S=\min\nolimits_{1\leq i\leq n}\{\mathbb{D}_i\}$, $T= \mathop
{\min}\nolimits_{1\leq i\neq j\leq n} \{ {d_{ij}\mathbb{D}_j^\alpha
\mathord{\left/
 {\vphantom {d_{ij}\mathbb{D}_j^\alpha \mathbb{D}_i^\alpha}} \right.
 \kern-\nulldelimiterspace} \mathbb{D}_i^\alpha}\}$.\\
\\
\textbf{Acknowledgements} This work was partially supported by the
National Natural Science Foundation of China (No. 11271334), the
Natural Science Foundation of Zhejiang Province, China (No.
LY12A01006).

\end{document}